\documentclass{dimacs-l}
\usepackage{amscd,amssymb,amsxtra}

\theoremstyle{plain}
\newtheorem{Theorem}{Theorem}[section]
\newtheorem{Observation}[Theorem]{Observation}

\newtheorem{Lemma}[Theorem]{Lemma}
\newtheorem{Sublemma}[Theorem]{Sublemma}
\newtheorem{Fact}[Theorem]{Fact}
\newtheorem{Corollary}[Theorem]{Corollary}
\newtheorem{Claim}[Theorem]{Claim}

\theoremstyle{definition}
\newtheorem{Definition}[Theorem]{Definition}
\newtheorem{Question}[Theorem]{Question}

\theoremstyle{remark}

\numberwithin{equation}{section}

\newcommand{\cf}{\operatorname{cf}}

\newcommand{\cal}{\mathcal}

\newcommand{\eop}{\bigstar}

\newcommand{\card}[1]{{\vert #1 \vert} }

\newenvironment{Proof}{\noindent{\bf Proof.}}{\par\bigskip} 

\newenvironment{Note1}{\noindent{\bf Note.}}{\par\bigskip} 

\newenvironment{Proof of the Subfact}
{\noindent{\bf Proof of the Subfact.}}{\par\bigskip}

\newenvironment{Proof of the Theorem}
{\noindent{\bf Proof of the Theorem.}}{\par\bigskip}

\newenvironment{Proof of Theorem 1}
{\noindent{\bf Proof of Theorem }}{\par\bigskip}

\newenvironment{Proof of the Conclusion}
{\noindent{\bf Proof of the Conclusion.}}{\par\bigskip}

\newenvironment{Proof of the Observation}
{\noindent{\bf Proof of the Observation.}}{\par\bigskip}

\newenvironment{Proof of the Fact}
{\noindent{\bf Proof of the Fact.}}{\par\bigskip}

\newenvironment{Proof of the Lemma}
{\noindent{\bf Proof of the Lemma.}}{\par\bigskip}

\newenvironment{Proof of the Claim}
{\noindent{\bf Proof of the Claim.}}{\par\bigskip}

\newenvironment{Proof of the Subclaim}
{\noindent{\bf Proof of the Subclaim.}}{\par\medskip}

\newenvironment{Proof of the Sublemma}
{\noindent{\bf Proof of the Sublemma.}}{\par\medskip}

\newenvironment{Proof of the Main Claim}
{\noindent{\bf Proof of the Main Claim.}}{\par\bigskip}

\newcommand{\elementary}{\prec}

\newcommand{\Bbf}{\Bbb}

\newcommand{\into}{\rightarrow}

\newcommand{\rest}{\upharpoonright}  

\newcommand{\satisfies}{\vDash}
\newcommand{\deq}{\buildrel{\rm def}\over =}

\newcommand{\BB}{{\cal B}}

\newcommand{\FF}{{\cal F}}
\newcommand{\GG}{{\cal G}}
\newcommand{\HH}{{\cal H}}

\newcommand{\PP}{{\cal P}}

\begin{document}

\title[On $D$-spaces and Discrete Families]
{On $D$-spaces and Discrete Families of Sets}

\author{Mirna D\v zamonja}
\address
{School of Mathematics\\
University of East Anglia\\
Norwich, NR4 7TJ,UK}
\email{M.Dzamonja@uea.ac.uk}

\thanks{The author thanks EPSRC, Leverhulme Trust and NATO for their
support through research grants. Warm thanks to William Fleissner for his many
useful comments, corrections and additions to the paper.}

\subjclass{03E35, 54E20,
03E55.}
\keywords{$D$-spaces, reflection, discrete families}


\begin{abstract}
We prove several reflection theorems on $D$-spaces, which are Hausdorff
topological spaces $X$ in which for every open neighbourhood assignment
$U$ there is a closed discrete subspace $D$ such that 
\[
\bigcup\{U(x):\,x\in D\}=X.
\]
The upwards reflection theorems are obtained in the presence of a forcing
axiom,
while most of
the downwards reflection results use large cardinal assumptions.

The combinatorial content of arguments showing that a
given space is a $D$-space, can be formulated using the concept of
discrete families. We note the connection between non-reflection arguments
involving discrete families and the well known question of the existence of
families allowing partial transversals without having a transversal themselves,
and use it to give non-trivial instances of the incompactness phenomenon
in the context of discretisations.
\end{abstract}

\maketitle

\section{Introduction} We prove some reflection
results about $D$-spaces, and note their combinatorial equivalent entitled
discrete families. $D$-spaces
were defined by E.K. van Douwen in 1978, and studied
by van Douwen and W.K. Pfeffer \cite{DoPf}, van Douwen and D. Lutzer
\cite{Lutzer} and W. Fleissner and A.M. Stanley \cite{Adrienne}, among others.
To define $D$-spaces, recall that an open neighbourhood assignment (ONA) in a topological
space $X$ is a function $U$ on $X$ such that for all $x\in X$ we have that
$U(x)$ is an open neighbourhood of $x$. A space is said to be a $D$-space iff
for every ONA $U$ of $X$, there is a closed discrete $D\subseteq X$
such that $\bigcup\{U(x):\,x\in D\}=X$. We equivalently say that a space is $D$ or
{\em has property} $D$.

It is easy to see that compact spaces, discrete spaces and metric spaces are
$D$. A very puzzling open question about $D$-spaces is if all Lindel\"of
spaces are $D$. In an attempt to solve this question van Douwen and Pfeffer
\cite{DoPf} studied the Sorgenfrey line $S$ and proved that all finite powers of
$S$ are $D$, as well as introducing a larger class of spaces which are $D$.
Continuing this line of research,
the known results about $D$-spaces often concentrate on
generalised metric spaces and
linearly ordered
spaces (LOTS) and their products. Much of this effort, including the question
about Lindel\"{o}f spaces being $D$, can be viewed as focussing on the natural
generalisations of the fact that compact spaces are $D$. We note that there is
another natural line of generalisation of the basic facts about $D$-spaces,
namely the observation that discrete spaces are $D$. In this vein, one
should consider spaces which are locally small, that is, every point has
a small neigbourhood.
We prove
for example
(Corollary \ref{xxx}) that it is consistent with $CH$ and
$2^{\aleph_1}>\aleph_2$ that every  locally countable Hausdorff space of size $\le\aleph_2$ in which every open subspace of size
$\le\aleph_1$ is $D$ in a strong sense, is $D$ itself. It is perhaps worth
mentioning that there seems to be inherent difficulties in proving consistency
results about $D$-spaces. In particular, the result we prove is to our
knowledge the first such result, and it is still a reflection argument rather
than an outright consistency result.

As the result just mentioned is an upward reflection result, it is natural to
ask if there are also downward reflection results. With the help of large cardinal
assumptions, we can get some such results. Namely, we prove that if $\kappa$ is
a measurable cardinal
and $\sigma<\kappa$, then every locally $<\kappa$ $\sigma$-$D$ space of size $\kappa$
in which every point has a point-base of size $<\kappa$, has open 
$\sigma$-$D$ subspaces
of sizes arbitrarily close to but less than $\kappa$
(the notion of $\sigma$-$D$ spaces is defined below and is crucial
for the upwards reflection results from the first section).

The same argument can be used with $D$ in place of $\sigma$-$D$,
but we give instead an improvement due to W. Fleissner, where
the large cardinal assumption is reduced to a strong inaccessible
and the assumption of small character is not needed. The point in
both arguments is to get
open subspaces, as getting closed subspaces, for instance, is very easy since
the property of being $D$ reflects downward to closed subspaces.

The paper finishes with a section on discrete families. When working with
$D$-spaces, one quickly realises that there is a combinatorial argument
repeatedly being used. Formalising the ingredients of this argument,
we can abstract the combinatorial content of the $D$-space context, and
arrive at the notion of discrete families. We give a short discussion of these,
and note that reflection arguments about discrete families have a lot to do
with the well studied problem of the existence of transversals. Then we
use this observation and the known results about transversals
to give non-trivial instances of the existence of discrete families.

Although many of the results mentioned here are still valid if
we work with spaces which are only assumed to be $T_1$, we shall for
simplicity only study Hausdorff
topological spaces.
\section{A consistency result on upwards reflection} In this section we
prove the upwards reflection theorem
announced in the introduction (Corollary \ref{xxx}), obtaining it as a
consequence of the following more general Theorem.

\begin{Theorem}\label{upwards} Suppose that ${\bf V}$ is a universe in which
\[
\aleph_0\le\lambda=\lambda^{<\lambda}<\lambda^+=\kappa<\lambda^{++}<\mu
\]
and $\mu^{\lambda^+}=\mu$, while $2^\lambda=\kappa$ and $2^{\kappa}=
\kappa^+$.

Then there is a cofinality and cardinality preserving forcing extension of
${\bf V}$ in which no bounded subsets are added to $\kappa$ and the following
hold: \begin{description}
\item{(i)} $2^\lambda=\lambda^+$ and $2^{\lambda^+}=\mu$,
\item{(ii)} every locally $<\kappa$ topological space of size $\le\kappa^+$
in which for every ONA there is a finer ONA with respect to
which all open subspaces
of size $\le\kappa$ are $\kappa$-$D$, is
$D$ itself.
\end{description}
\end{Theorem}
\rightline{$\between$}

Let us first recall the definitions of an ONA and a
$D$-space,
which were mentioned in the Introduction, give some background to the
concepts needed for the proof, and most importantly, define what
a $\kappa$-$D$-space is.

\begin{Definition}\label{ona}(1)
An {\em open neighbourhood assignment} (ONA) in a topological
space $(X, \tau)$ is a function $U:\,X\into\tau$ such that for all $x\in X$ we
have that $x\in U(x)$.

{\noindent(2)} A space is said to be a $D$-{\em space} iff
for every ONA $U$ of $X$, there is a closed discrete $D\subseteq X$
such that $\bigcup\{U(x):\,x\in D\}=X$.

{\noindent(3)} For a cardinal $\kappa$, a topological space $X$ is said to be
{\em locally $<\kappa$} iff there is an ONA $U$ of $X$ such that $\card{U(x)}
<\kappa$ for all $x\in X$.

{\noindent(4)} If $U$ is an ONA of $X$ for which there is a closed discrete $D$
with \[
\bigcup U``D\deq\bigcup\{U(x):\,x\in D\}=X,
\]
we say that $X$ is $D$ with respect to $U$.
\end{Definition}

\begin{Observation}\label{refinement} If $U, U^\ast$ are ONA of $X$ such that
$\forall x(U(x)\subseteq U^\ast(x))$, and $X$ is $D$ with respect to $U$, then
$X$ is $D$ with respect to $U^\ast$.
\end{Observation}
\rightline{$\between$}

Given a space $X$ and an ONA $U$ on it, if one tries to construct inductively
or otherwise a subspace $D\subseteq X$ demonstrating that $X$ is a $D$-space
with respect to $U$,
there are two apparent difficulties that one may run into. One of these is that
taking unions of infinitely many closed discrete subspaces does not necessarily give a
closed subspace. This difficulty is resolved through the use of $U$-sticky sets,
as introduced and studied by
Fleissner and Stanley in \cite{Adrienne}.

\begin{Definition}\label{PU} Given an
ONA $U$ of a topological
space $X$
\begin{description}
\item{(1)} A subspace $D$ of $X$ is said to be $U$-{\em sticky} iff $D$ is
closed discrete and satisfies
\[
(\forall x\in X) [U(x)\cap D\neq\emptyset\implies x \in \bigcup U``D].
\]
\item{(2)} The partial order ${\Bbf P}_U={\Bbf P}_U(X)$ is defined by
letting
\[
{\Bbf P}_U\deq\{D\subseteq X:\,D\mbox{ is }U\mbox{-sticky}\},
\]
ordered by letting $D\le D'$ (where $D'$ is a stronger condition) iff
$D\subseteq D'$ and $(D'\setminus D)\cap \bigcup U``D=\emptyset$.
\end{description}
\end{Definition}
\rightline{$\between$}

Fleissner and Stanley proved that ${\Bbf P}_U$ is well behaved with 
respect to the unions of $\le$-increasing chains, see the following Theorem
\ref{FlSt}(1). This gives hope that one could use ${\Bbf P}_U$ in inductive
constructions, or as a forcing notion, but at least as much as the latter is
concerned, this hope is slighted by a further result of Fleissner and Stanley.
Namely, the second of the difficulties mentioned above, is a
density problem: given a $U$-sticky $D$, and an $x\in X$, can we find an
extension $D'$ of $D$ within ${\Bbf P}_U$ for which we have
$x\in \bigcup U``D'$? The second part of 
Theorem \ref{FlSt} shows that such a density condition is present iff $X$ is
a $D$-space. While this gives a very interesting characterisation of
$D$-spaces, it also shows that ${\Bbf P}_U(X)$ cannot be used as a forcing
notion to make $X$ into a $D$-space, if $X$ was not a $D$-space to start with.

\begin{Theorem}\label{FlSt} [Fleissner-Stanley] Given a topological
space $X$ and an open neighbourhood assignment $U$ of $X$. Let ${\Bbf P}_U(X)$
be as defined above, in Definition \ref{PU}(2). Then
\begin{description}
\item{(1)} if ${\cal D}$ is a subset of ${\Bbf P}_U$ in which every $D, D'$
satisfy either $D\le D'$ or $D'\le D$, then $\bigcup {\cal D}\in {\Bbf P}_U$.
\item{(2)} $X$ is a $D$-space iff for every ONA $U$ of $X$, every $D\in {\Bbf
P}_U$,
and every $x\in X$, there is a $D'\ge D$ in ${\Bbf P}_U$ such that $x\in 
\bigcup U``D'$.
\end{description}
\end{Theorem}
\rightline{$\between$}

In the following discussion we shall be assuming that $X$ is a given
topological
space and $U$ an ONA on $X$.
We shall work with a variant of $U$-stickiness which will be used as a forcing
notion.  For a given regular cardinal $\kappa\ge \aleph_1$ we define the
partial order  ${\Bbf P}_U^\kappa$ as follows:

\begin{Definition}\label{PUkappa} The partial order
${\Bbf P}_U^\kappa={\Bbf
P}_U^\kappa(X)$ is defined by letting
\[
{\Bbf P}_U^\kappa\deq\{D\subseteq X:\,D\mbox{ is
}U\mbox{-sticky}\,\,\&\,\,\card{D}<\kappa\},
\]
ordered by letting $D\le D'$ (where $D'$ is a stronger condition) iff
$D\subseteq D'$ and $(D'\setminus D)\cap \bigcup U``D=\emptyset$.
\end{Definition}
\rightline{$\between$}
We start the discussion of ${\Bbf P}_U^\kappa$ by a slight generalisation of
Theorem \ref{FlSt}(1), proved in a manner similar to the one used for the
proof of that theorem.

\begin{Observation}\label{easy} (1) If ${\cal D}$ is a
subset of ${\Bbf P}_U^\kappa$  with $\card{{\cal D}}<\kappa$
and such that for each
$D,D'\in {\cal D}$ there is a $D''\in {\cal D}$ with $D''\ge D,D'$,
then $\bigcup {\cal D}$ is
an element of ${\Bbf P}_U^\kappa$.

{\noindent (2)} If ${\cal D}$ is a directed subset of ${\Bbf P}_U$, then
$\bigcup {\cal D}$ is an element of ${\Bbf P}_U$.
\end{Observation}

\begin{Proof} (1) Let $D^\ast\deq\bigcup{\cal D}$, and 
let $x\in X$. If $U(x)\cap D^\ast=\emptyset$, then
certainly $x\notin \overline{D^\ast\setminus\{x\}}$. Otherwise, there is $D\in
{\cal D}$ such that $U(x)\cap D\neq\emptyset$. Since $D\in {\Bbf P}_U^\kappa$,
we have that $x\in \bigcup U``D\subseteq \bigcup U``D^\ast$. This demonstrates the
second part of the definition of being $U$-sticky. 

Fix an $x\in X$ again. We shall use it to show
that $D^\ast$ is closed and discrete. Suppose that $U(x)\cap D\neq\emptyset$
for some $D\in {\cal D}$. 
Since $D$ is closed discrete, there is an open
neigbourhood $V$ of $x$ with $V\cap D$ finite. As $U(x)\cap D\neq\emptyset$
we have that $x\in \bigcup U``D$, so we can assume that $V\subseteq U(x)\cap
\bigcup U``D$. Given $D'\in {\cal D}$, let $D''$ be a common extension of $D$
and $D'$ in  ${\Bbf P}_U^\kappa$. Hence
\[
(D'\setminus D)\cap \bigcup U``D\subseteq (D''\setminus D)\cap \bigcup U``D=\emptyset,
\]
so $V\cap D'\subseteq \bigcup U``D\cap D'\subseteq \bigcup U``D\cap D$,
and so $V\cap D'\subseteq V\cap D$. In
conclusion, $V\cap D^\ast=V\cap D$ is finite, and so $x\notin
\overline{D^\ast\setminus\{x\}}$, by the Hausdorff property of $X$. If
$U(x)\cap D=\emptyset$ for all $D\in {\cal D}$, then $U(x)\cap
D^\ast=\emptyset$, so clearly $x\notin \overline{D^\ast\setminus\{x\}}$. This argument demonstrates
that $\bigcup{\cal D}$ is closed discrete. As $\card{{\cal D}}<\kappa$, we
have $\card{D^\ast}<\kappa$.

{\noindent (2)} The same argument as above, omitting the last sentence.
$\eop_{\ref{easy}}$
\end{Proof}

The order ${\Bbf P}_U^\kappa$ can be used to define what is meant by a
$\kappa$-$D$-space, keeping until a further notice the convention that $\kappa$ is
a chosen uncountable regular cardinal.

\begin{Definition}\label{kappaD} (1) $X$ is said to be $\kappa$-$D$
{\em with respect to an} ONA $U$ {\em of} $X$
iff for every $D\in {\Bbf P}_U^\kappa$, and
every $x\in X$, there is $D'\ge D$ in ${\Bbf P}_U^\kappa$ such that $x\in \bigcup
 U``D'$.

{\noindent (2)} Given two ONA $U$ and $V$ of $X$, we say that $V$ is
{\em finer than} $U$ iff $V(x)\subseteq U(x)$ for all $x\in X$.

{\noindent (3)} $X$ is $\kappa$-$D$ iff for every ONA $U$ of $X$ there is a
finer ONA $V$ such that $X$ is $\kappa$-$D$ with respect to $V$.

We use the terminology ``strongly-$D$" in place of ``$\aleph_1$-$D$".
\end{Definition}

\begin{Note1} The definition of $X$ being $\kappa$-$D$ does not require
$X$ to be $\kappa$-$D$ with respect to every ONA $U$ of $X$. Note also that $X$
is $D$ iff $X$ is $\card{X}^+$-$D$.
\end{Note1}
\rightline{$\between$}

The choice of our terminology can be explained by the following

\begin{Observation}\label{explanation} Suppose that
$\kappa\ge\aleph_1$ is regular. Then every $\kappa$-$D$-space of
size $\le \kappa$ is $D$. 
\end{Observation}

\begin{Proof} Suppose
that $X=\{x_\alpha:\,\alpha<\alpha^\ast\le\kappa\}$ is a
given $\kappa$-$D$-space
and that $U^\ast$ is a given ONA of $X$. Let $U$ be a finer ONA such that
$X$ is $\kappa$-$D$ with respect to $U$. By
induction on $\alpha$ we define $\langle D_\alpha:\,\alpha\le\alpha^\ast\rangle$, a continuous increasing
chain of elements of ${\Bbf P}_U^\kappa$, with $D_0=\emptyset$ and $x_\alpha\in
\bigcup U``D_{\alpha+1}$, while $\bigcup_{\alpha\le\alpha^\ast}U``D_\alpha=X$. The
induction at successor stages uses the assumption of $\kappa$-$D$-ness, and at limit
stages, Observation \ref{easy}(1). Using the same Observation, we can see that
taking
$D\deq\bigcup_{\alpha\le\alpha^\ast}D_\alpha$ demonstrates that $X$ is $D$ with
 respect  to  $U$, hence $X$ is $D$  with
 respect  to  $U^\ast$ (Observation \ref{refinement}).
 $\eop_{\ref{explanation}}$ \end{Proof}

Further discussion about the relationship between being $D$ and $\kappa$-$D$
can be found at the end of this section.
We intend to use ${\Bbf P}_U^\kappa$ in a universe of set theory in which a
certain version of Martin's axiom for $\kappa^+$ holds. For this we shall need a
theorem of S. Shelah from \cite{Sh80} (there proved with $\kappa=\aleph_1$,
but, as is well known to the author of \cite{Sh80}  and has been used in many of his
results, the same proof gives the more general result here quoted as Theorem
\ref{Shelah}). To
introduce this theorem, we need a definition.

\begin{Definition}\label{star} We say that a forcing notion ${\Bbf P}$
satisfies $(\ast_\kappa)$ iff the following conditions (a)-(c) hold:
\begin{description}
\item{(a)} if $p,q$ are compatible in ${\Bbf P}$, then they have a least upper
bound (lub) in ${\Bbf P}$,
\item{(b)} if $\langle p_\alpha:\,\alpha<\alpha_\ast<\kappa\rangle$ is an
increasing sequence of conditions in ${\Bbf P}$, then the sequence
has the lub in ${\Bbf P}$,
\item{(c)} if $\{p_i:\,i<\kappa^+\}$ is a set of conditions in ${\Bbf P}$,
then there is a club $C\subseteq\kappa^+$ and a regressive function
$f$ on $C$ such that whenever $i,j\in C$ are of cofinality $\kappa$ and
$f(i)=f(j)$, then $p_i$ and $p_j$ are compatible.
\end{description}
\end{Definition} 

\begin{Theorem}\label{Shelah} [Shelah] Suppose that ${\bf V}$ is a universe 
of set theory in which
\[
\aleph_0\le\lambda=\lambda^{<\lambda}<\lambda^+=\kappa<\lambda^{++}<\mu
\]
and $\mu^{\lambda^+}=\mu$, while $2^\lambda=\kappa$ and $2^{\kappa}=
\kappa^+$.

Then there is a cofinality and cardinality preserving forcing extension of
${\bf V}$ in which no bounded subsets are added to $\kappa$ and the following
hold: \begin{description}
\item{(i)} $2^\lambda=\lambda^+$ and $2^{\kappa}=\mu$,
\item{(ii)} for every forcing ${\Bbf P}$ satisfying $(\ast_\kappa)$ and
a collection ${\cal D}=\{D_\zeta:\,\zeta<\zeta^\ast<\mu\}$ of dense subsets of
${\Bbf P}$, there is a filter $G$ of ${\Bbf P}$ with $G\cap
D_\zeta\neq\emptyset$ for every $\zeta$.
\end{description}
\end{Theorem}
\rightline{$\between$}

Observation \ref{easy} allows for an easy proof of the following

\begin{Observation}\label{easy2} Suppose that $\langle D_\alpha:\,\alpha<
\alpha^\ast< \kappa\rangle$ is an increasing sequence of conditions in ${\Bbf P}^\kappa_U$.
Then $\bigcup_{\alpha<\alpha^\ast} D_\alpha$ is the lub of the sequence.
\end{Observation}

\begin{Proof} By Observation \ref{easy}(1) we have $D\deq
\bigcup_{\alpha<\alpha^\ast} D_\alpha\in {\Bbf P}^\kappa_U$, so we only need to show that 
$D$ is an extension
of each $D_\alpha$, and that is actually the least such condition. We can without loss of generality
assume that $\alpha^\ast$ is a limit ordinal.

Let $\alpha<\alpha^\ast$. If $x\in D\setminus D_\alpha$, then for some $\beta\in
(\alpha,\alpha^\ast)$ we have $x\in D_\beta$, so, since $D_\alpha\le D_\beta$, we have that
$x\notin \bigcup U``D_\alpha$. Hence $D_\alpha\le D$.

Suppose that $D'\ge D_\alpha$ for all $\alpha$. In particular, $D'\supseteq
D$. If $x\in D'\setminus D$ but $x\in \bigcup U``D$, then there is
$\alpha<\alpha^\ast$ such that $x\in \bigcup U``D_\alpha$, contradicting that
$D_\alpha\subseteq D$ and $D_\alpha\le D'$. Hence $x\in D'\setminus D\implies
x\notin \bigcup U``D$, so $D\le D'$.
$\eop_{\ref{easy2}}$
\end{Proof}

The next item needed for the proof of Theorem \ref{upwards} is the existence of
a lub of two conditions compatible in ${\Bbf P}_U^\kappa$, and the proof of
this is another elementary argument of the sort used to prove the previous
Observations.

\begin{Observation}\label{easy1} If $D,D'\in {\Bbf P}_U^\kappa$, then $D$ and
$D'$ are compatible iff
their union $D\cup D'$ is their common upper bound, in which case
$D\cup D'$ is the lub of $D$ and $D'$ in ${\Bbf P}_U^\kappa$.
\end{Observation}

\begin{Proof} Let us prove the nontrivial part of this Observation, so assume
that $D''$ is a common upper bound of $D$ and $D'$. Clearly, $D\cup D'$ is
closed and discrete, and has size $<\kappa$. If $x$ is an element of
$(D\cup D')\setminus D$, then $x\in D''\setminus D$, and in particular
$x\notin \bigcup U``D$. By symmetry, the same argument can be applied with $D'$ in
place of $D$. Next, if $x$ is such that $U(x)\cap (D\cup D')\neq\emptyset$,
then either $U(x)\cap D\neq\emptyset$, hence $x\in \bigcup U``D$, or similarly,
$x\in \bigcup U``D'$.

Finally, we have to show that $D''\ge D\cup D'$, which can be done by similar
elementary arguments.
$\eop_{\ref{easy1}}$
\end{Proof}

Now we are ready to prove

\begin{Lemma}\label{cc} Let $X$ be a topological space whose points
are among the ordinals $\le\kappa^+$.
Suppose that $U$ is an ONA of $X$ that has the
property $\card{U(x)}<\kappa \,(\forall x\in X)$, and assume that 
$\kappa^{<\kappa}=\kappa$. Then
${\Bbf P}^\kappa_U$ satisfies $(\ast_\kappa)$.
\end{Lemma}

\begin{Proof} Observations \ref{easy1} and \ref{easy2} provide us with the
properties (a) and (b) from Definition of $(\ast_\kappa)$. We shall now prove
the required chain condition. Suppose that we are given
$\{D_\alpha:\,\alpha<\kappa^+\}$ from ${\Bbf P}^\kappa_U$. As
$\kappa^{<\kappa}=\kappa$, we also have
$(\kappa^+)^{<\kappa}=\kappa^+$, so we can fix a bijection $F$ from 
$([\kappa^+]^{<\kappa})^2$ onto $\kappa^+$. We now define several subsets of
$\kappa^+$:
\[
C_0\deq\{\alpha<\kappa^+:\,(\forall\beta<\alpha)(\forall A,B\in
[\beta]^{<\kappa})\,F(A,B)<\alpha\},
\]
\[
C_1\deq\{\alpha<\kappa^+:\,(\forall \beta<\alpha)D_\beta\subseteq\alpha\}
\]
and
\[
C_2\deq
\{\alpha<\kappa^+:\,(\forall\beta<\alpha) U(\beta)\subseteq\alpha\}.
\]
Let $C\deq (C_0\cap C_1\cap C_2)\setminus \{0\}$. Standard arguments
show that $C$ is a club of
$\kappa^+$. This will be the club demonstrating the required condition.
In order to finish the demonstration, we also need to define a 
regressive function $f$. To motivate its definition, let us first prove

\begin{Sublemma}\label{sublemma} Suppose $\alpha<\beta<\kappa^+$ are such that
\[
D_\beta\cap \bigcup U``(\bigcup_{\gamma<\beta} D_\gamma)=
D_\alpha\cap \bigcup U``(\bigcup_{\gamma<\alpha} D_\gamma)
\]
and
\[
(\bigcup U``D_\beta)\cap \bigcup_{\gamma<\beta} D_\gamma=
(\bigcup U``D_\alpha)\cap \bigcup_{\gamma<\alpha} D_\gamma.
\]
Then $D_\alpha$ and $D_\beta$ are compatible.
\end{Sublemma}

\begin{Proof of the Sublemma} Let $\alpha$ and $\beta$ be as claimed.
We shall show that $D\deq D_\alpha\cup D_\beta$ is a common upper bound
of $D_\alpha$ and $D_\beta$.
$D$ is clearly a closed and discrete
superset of $D_\alpha$ and $D_\beta$, and has size $<\kappa$. Suppose that
$x\in X$ and $U(x)\cap D\neq\emptyset$, then $U(x)\cap D_l\neq\emptyset$
for some $l\in \{\alpha,\beta\}$. In any case, $x\in \bigcup U``D_l
\subseteq \bigcup U``D$.

If $x\in D\setminus D_\alpha$, then $x\in D_\beta\setminus D_\alpha$. Should
$x$ belong to $(D_\beta\setminus D_\alpha)\cap \bigcup U``D_\alpha$,
we would have 
\[
x\in D_\beta\cap \bigcup U``(\bigcup_{\gamma<\beta}D_\gamma)=
D_\alpha\cap \bigcup U``(\bigcup_{\gamma<\alpha} D_\gamma),
\]
a contradiction, since $x$ is assumed not to be in $D_\alpha$. Hence,
the intersection between
$(D\setminus D_\alpha)$ and $\bigcup U``D_\alpha$ is empty. As we have
assumed $\alpha<\beta$, this argument does not automatically yield the
analogous conclusion with $\beta$ in place of $\alpha$, but the rest of our
assumptions about $\alpha$ and $\beta$ can be used now, as follows.

If $x\in D\setminus D_\beta$, then $x\in D_\alpha\setminus D_\beta$. Supposing
that also $x\in \bigcup U``D_\beta$, we have $x\in \bigcup U``D_\beta\cap
\bigcup_{\gamma <\beta}D_\gamma$, which is the same as $\bigcup U``D_\alpha\cap
\bigcup_{\gamma<\alpha} D_\gamma$. Hence 
\[
x\in D_\alpha\cap\bigcup_{\gamma
<\alpha}D_\gamma\subseteq D_\alpha\cap \bigcup
U``(\bigcup_{\gamma<\alpha}D_\gamma)=D_\beta\cap\bigcup U``(\bigcup_{\gamma<\beta}
D_\beta),
\]
contradicting the assumption that $x\notin D_\beta$.
$\eop_{\ref{sublemma}}$
\end{Proof of the Sublemma}

{\noindent{\em Proof of Lemma \ref{cc} continued.}
We define $f$ on $C$ by letting for $\alpha\in C$ with 
$\cf(\alpha)=\kappa$
\[
f(\alpha)\deq\ F(D_\alpha\cap\bigcup U``(\bigcup_{\gamma<\alpha}D_\gamma),
\bigcup U``D_\alpha\cap\bigcup_{\gamma<\alpha}D_\gamma),
\]
and letting $f(\alpha)=0$ otherwise. Notice that the cardinal assumptions
on $U$ and $D_\alpha$s guarantee that $f$ is well defined.
Let us show that it is regressive on $C$. For the nontrivial part of this,
let $\alpha\in C$ be of cofinality $\kappa$. By the choice of $C_1$ and $C_2$
we have that
\[
\bigcup_{\gamma<\alpha}
D_\gamma\subseteq\alpha
\mbox{ and }\bigcup U``(\bigcup_{\gamma<\alpha} D_\gamma) \subseteq\alpha.
\]
Hence $A\deq D_\alpha\cap [\bigcup U``(\bigcup_{\gamma<\alpha} D_\gamma)]$
and $B\deq (\bigcup_{\gamma<\alpha}
D_\gamma)\cap \bigcup U`` D_\alpha$ are both subsets of
$\alpha$ of size $<\kappa$, hence bounded. By the choice of $C_0$ we have
$f(\alpha)=F(A,B)<\alpha$.

Now let us see that $C$ and $f$ work: suppose that $\alpha<\beta$ are in $C$,
have cofinality $\kappa$ and $f(\alpha)=f(\beta)$. By Sublemma
\ref{sublemma} we have that $D_\alpha$ and $D_\beta$ are compatible.
$\eop_{\ref{cc}}$
}\end{Proof}

\begin{Proof of Theorem 1}{\bf \ref{upwards}.} Starting with ${\bf V}$
and the cardinals as in the statement of the
Theorem, using Shelah's Theorem we pass to a universe $W$ in which the 
conclusions of that theorem hold. From now on, let us work in $W$. Let
$(X,\tau)$ be a given locally $<\kappa$-space of
size $\le\kappa^+$, 
such that
for every ONA $U^\ast$ of $X$ there is a finer ONA $U$ with the
property that every open
subspace of size $\le\kappa$ of $X$ is $\kappa$-$D$ with respect to $U$. We
shall show that $X$ is a $D$-space. By Observation \ref{explanation}, we can assume that the size of $X$ is
$\kappa^+$. Without loss of generality the points of $X$ are the ordinals $<\kappa^+$.

Let $\{V(x):\,x\in X\}$ be an open neigbourhood assignment on $X$ which
demonstrates that $X$ is locally $<\kappa$, hence $\card{V(x)}<\kappa$
for every $x\in X$. Note that in order to show that $X$ is $D$, we may
concentrate on those ONA $U^\ast$ of $X$ for which we have $U^\ast(x)\subseteq
V(x)$ for all $x\in X$. Let $U^\ast$ be such an ONA
and let $U$ be a finer ONA with respect to which all open subspaces of
size $\le\kappa$ are $\kappa$-$D$. It suffices to show that $X$ is $D$
with respect to $U$.

Let
${\Bbf P}={\Bbf P}^\kappa_U$. As we have assumed that $\kappa^{<\kappa}=\kappa$ in ${\bf V}$, and no bounded
subsets of $\kappa$ are added when $W$ is formed, we have that 
$\kappa^{<\kappa}=\kappa$ holds in $W$. Therefore, Lemma \ref{cc} applies,
and we conclude that ${\Bbf P}$ satisfies $(\ast_\kappa)$. The main part of the
rest of the proof is a density argument.

\begin{Claim}\label{main} For every $x\in X$, the set
\[
{\cal D}_x\deq\{D\in {\Bbf P}:\,x\in \bigcup U``D\}
\]
is dense in ${\Bbf P}$.
\end{Claim}

\begin{Note1} One may wonder if the fact that Claim \ref{main} holds and
the part (2) of Fleissner-Stanley Theorem, do not automatically imply that $X$
is $D$, without a reference to forcing. This is not necessarily the case, as
Claim \ref{main} only refers to $U$-sticky sets of size $<\kappa$.
\end{Note1}

\begin{Proof of the Claim} 
Let us fix a large enough regular cardinal $\chi$.

Given $D\in {\Bbf P}$ and $x\in X$. We choose $N\elementary (\HH(\chi),
\in, <^\ast_\chi)$ of size $\kappa$ and such that
\begin{description}
\item{(i)} $X,\tau, D, x, U \in N$,
\item{(ii)} ${}^{\kappa>}N\subseteq N$ and $\kappa+1\subseteq N$.
\end{description}
Such a choice is possible, as $\kappa^{<\kappa}=\kappa$.
We shall look for a $D'$ which is a required extension of $D$ in ${\cal D}_x$.

Consider first the subset $N\cap X$ of $X$ in the subspace topology. If $y\in
N\cap X$, then $U(y)\in N$. As $\HH(\chi)\models``\card{U(y)}<\kappa"$,
there is an $\alpha<\kappa$ and a function $f$ from $\alpha$ onto
$U(y)$. By elementarity, there is such a function in $N$, and since
$\alpha\subseteq N$, we have that $U(y)\subseteq N$. This shows that
$N\cap X$ is an open subspace of $X$, and hence, by our assumptions,
$N\cap X$ is $\kappa$-$D$ with respect to $U$.
In order to use this fact, we shall show that $D\in 
{{\Bbf P}}^\kappa_{U^N\rest(N\cap X)} (N\cap X)$.
Here $U^N\rest(N\cap X)$ stands for the function assigning $U(x)\cap N
(=U(x))$
to each $x\in N\cap X$. Note
again that we are considering $N\cap X$ in the subspace topology,
not the
topology induced by elementarity, so in particular $U^N\rest(N\cap X)$ is an
ONA of $N\cap X$. Let ${\Bbf R}\deq {\Bbf P}^\kappa_{U^N\rest(N\cap X)} (N\cap
X)$.

An argument similar to the one showing that $N\cap X$ is open, shows that
$D$ is a subset of $N \cap X$. Clearly, $D$ is closed and discrete in $N\cap X$,
and has size $<\kappa$. Suppose that $y\in N\cap X$ is such that $U(y)\cap
D\cap (N\cap X)
\neq\emptyset$, then by the fact that $D\in{\Bbf P}$, we have that $y\in
\bigcup U``D$, so $y\in \bigcup \{U(z)\cap N:\,z\in D\}$. This demonstrates that 
$D\in {\Bbf R}$. Hence, by our assumptions
(as $x\in N$), there is
$D'\in {\Bbf R}$ with $D\le_{{\Bbf R}}
 D'$ and $x\in \bigcup\{U(y)\cap N:\,y\in D'\}$.
In particular, $D\subseteq D'$ and $\card{D'}<\kappa$, and $x\in \bigcup
U``D'$. As $D'\in [N]^{<\kappa}$, we have that $D'\in N$.
It is easily seen that, since $N\cap X$ is open, the fact that $D'$ is
discrete in $N\cap X$, implies that $D'$ is discrete in $X$.

To show that $D'$ is closed in $X$, we shall have to use elementarity.
If $y\in (N\cap X)\setminus D'$, we have that for some open $O$ containing
$y$ we have $O\cap (N\cap X)\cap D'=\emptyset$, as $D'$ is closed in $N\cap X$.
By $X$ being locally $<\kappa$, we can assume that $\card{O}<\kappa$. Letting
$V=O\cap (N\cap X)$, we get that $y\in V\in \tau$ and $V$ is a subset of $N$
of size $<\kappa$, hence $V\in N$. Consequently,
\[
N\satisfies``(\forall y \in (X\setminus D'))(\exists V\in \tau)(y\in
V\,\,\&\,\,V\cap D'=\emptyset)",
\]
so the same is true in $\HH(\chi)$, demonstrating that $D'$ is
closed in $X$.

Note that if $y\in D'\setminus D$, then, using that $D'\subseteq N$
and $D\le_{\Bbf R} D'$, we have
that $y\notin \bigcup\{U(z)\cap N:\,z\in D\}$, so $y\notin \bigcup U``D$. 
For the rest of the proof, suppose that
$z\in X$ is such $U(z)\cap D'\neq\emptyset$, yet $z\notin \bigcup
U``D'$. By elementarity, there is $z'\in N$ such that $U(z')\cap
D'\neq\emptyset$ and $z'\notin \bigcup\{U(y)\cap N:\,y\in D'\cap N\}$. As
$D'\subseteq N$, we have $z'\notin \bigcup\{U(y)\cap N:\,y\in D'\}$,
contradicting the fact that $D'\in {\Bbf R}$. Hence $D'$ is as required.
$\eop_{\ref{main}}$
\end{Proof of the Claim}

{\noindent{\em Proof of Theorem \ref{upwards} continued.}}
By the choice of $W$, we can find a ${\Bbf P}$-filter $G$ such that
$G\cap{\cal D}_x\neq\emptyset$ for every $x\in X$. In particular, as $G$ is
a directed subset of ${\Bbf P}_U$,  by Observation \ref{easy}(2) we have that
$D\deq\bigcup G$ is an element of ${\Bbf P}_U$. Then $D$ is closed and
discrete, while the choice of the dense sets intersected by $G$ guarantees that
$\bigcup U``D=X$. $\eop_{\ref{upwards}}$
\end{Proof of Theorem 1}

\begin{Corollary}\label{xxx} It is consistent with $ZFC$ that $CH$ holds,
$2^{\aleph_1}>\aleph_2$ and every locally countable space of size $\le
\aleph_2$ in which
for every ONA $U^\ast$ there is a finer ONA with respect to which every open
subspace of size $\le\aleph_1$ is strongly $D$, is $D$ itself.
\end{Corollary}

\begin{Proof} We apply Theorem \ref{upwards}
with $\lambda=\aleph_0$.
$\eop_{\ref{xxx}}$
\end{Proof}

Having finished the proof of Theorem \ref{upwards}, there are several questions
that come to mind. Firstly, is there a difference between spaces which are
$\kappa$-$D$ and those which are simply $D$, and what does the assumption of
being locally $<\kappa$ contribute to this difference?
The simplest instance of this question would be:

\begin{Question}\label{prvo} Is there a locally countable $D$-space of
size $\aleph_1$ which is not strongly $D$?
\end{Question}

A simple argument shows that the simplest example of a locally countable
$D$-space of size $\aleph_1$ is strongly $D$, namely

\begin{Claim}\label{example} Suppose that $X$ is a non-stationary subset of
$\omega_1$, with  the order topology. Then $X$ is strongly $D$.
\end{Claim}

\begin{Note1} By the van Douwen-Lutzer \cite{Lutzer} characterisation of
linearly ordered $D$-spaces, such an $X$ is necessarily $D$.
\end{Note1}

\begin{Proof} Since $X$ is non-stationary, there is a club $C$ of $\omega_1$
with $C\cap X=\emptyset$. We can assume that $\min(C)<\min(X)$.
For
$\alpha\in X$ a limit ordinal, define $\beta_\alpha\deq\sup(C\cap\alpha)$, hence $\beta_\alpha<\alpha$ for every
such $\alpha$. Observe that,
since $C$ is unbounded, there is for any $\beta<\omega_1$ a
$\delta=\delta(\beta)$ such that for all $\alpha\ge\delta$, if $\beta_\alpha$
is defined, then $\beta_\alpha >\beta$.

Suppose that $U^\ast$ is an ONA of $X$. We choose a finer ONA such
that $U(\alpha)\subseteq (\beta_\alpha,\alpha+1)$ for limit ordinals
$\alpha\in X$, and $U(\alpha)=\{\alpha\}$ otherwise. Let $D\in
{\Bbf P}_U^{\aleph_1}(X)$ and $x\notin \bigcup U``D$. As $X$ is $D$, by
the Fleissner-Stanley Theorem, there is $D'\ge D$ in ${\Bbf P}_U(X)$,
with $x\in \bigcup U``D'$. Let $\gamma_0\in D'$ be such that $x\in U(\gamma_0)$,
and then define by induction an increasing sequence $\langle \gamma_n:\,n
<\omega\rangle$ of countable ordinals such that
$\gamma_{n+1}\ge\delta(\gamma_n)$, and $(\gamma_n,\gamma_{n+1})\cap
C\neq\emptyset$. Also require that $D\subseteq\gamma_1$. Now let
$\gamma\deq\sup_n\gamma_n$ and $D''=D'\cap \gamma$. As $\gamma\in C$, we have
that $\gamma\notin X$, hence $D''$ is closed in $X$.

Clearly, $D''$ is discrete, countable, satisfies $D\subseteq D''$ and $x\in 
\bigcup U``D''$. To show $D''\in {\Bbf P}^\kappa_U$, suppose $y\in X$ is such that
$U(y)\cap D''\neq\emptyset$. Then $y\in \bigcup U``D'$, by the choice of $D'$.
Suppose $y\notin \bigcup U``D''$, and let $\alpha\in D'$ be the minimal such that
$y\in U(\alpha)$. Hence $\alpha>\gamma$. If $\alpha$ is a successor, then
$y=\alpha$,
so $U(y)=\{\alpha\}$, contradicting the assumption $U(y)\cap D''\neq\emptyset$.
Hence, $\alpha$ is a limit and $y\in U(\alpha)\subseteq
(\beta_\alpha,\alpha+1)$.
As $\gamma\in C$, we have that $\beta_\alpha\ge\gamma$. Either $y$ is a
successor ordinal in $(\beta_\alpha,\alpha+1)$, contradicting $U(y)\cap
D''\neq\emptyset$, or $y$ is a limit ordinal. In the latter case, $\beta_y\ge
\beta_\alpha$, again contradicting $U(y)\cap
D''\neq\emptyset$.

Finally, $D\le D''$ because $D\le D'$.
$\eop_{\ref{example}}$
\end{Proof}

In fact, much more is true: by analysing the proof of Theorem 3.3. of
Fleissner-Stanley's paper \cite{Adrienne}, we can see that

\begin{Fact} A linearly ordered
topological space is $D$ iff it is $\kappa$-$D$ for all
regular uncountable $\kappa$.
\end{Fact}

W. Fleissner proved that the Cantor tree (see \cite{counter} for details) is
$D$ and not strongly $D$. The size of this space is $2^{\aleph_0}$.

Another question that might be worth asking is if the assumptions of Theorem
\ref{upwards}
are necessary. If one considers a non-reflecting stationary
subset $S$ of $\omega_2$ in the order topology, one has a space all of
whose subspaces of size $\le\aleph_1$ are strongly $D$, yet the space itself
is not $D$. If the points in $S$ have countable cofinality, this
space is even locally
countable. However, for a given ONA $U^\ast$ of
$X$ and a subspace $Y$ of $X$ with $\card{Y}\le\aleph_1$, the finer ONA $U$
with respect to which $Y$ is strongly $D$, depends on $Y$, i.e. cannot be
chosen uniformly for all $Y$, as in the assumptions of Theorem \ref{upwards}.
This indicates that some assumption additional to small
open subspaces being strongly $D$ is necessary in
the  statement of Corollary \ref{xxx}, and similarly in that of Theorem \ref{upwards}.

A tension between the existence of small open neighbourhoods and a certain
amount of compactness is a well studied subject, see for example
I. Juhasz, S. Shelah and L. Soukup's \cite{JuShSo}. Along these lines, one may
ask when there are locally countable non-discrete $D$-spaces of large
cardinality, although the fact that the relationship between being $D$
and other versions of compactness is not entirely clear, may mean that such a
question is premature. 

We end the section by discussing the possibility of strengthening Observation
\ref{explanation}.

\begin{Claim}\label{closedkappaD} Suppose that $X$ is a $\kappa$-$D$ space
and $Y$ is a closed subspace of $X$. Then $Y$ is $\kappa$-$D$.
\end{Claim}

\begin{Proof of the Claim} Let $\kappa, X, Y$ be as in the statement of the
Claim, and let $U^\ast$ be an ONA of $Y$. For $x\in X\setminus Y$, let
$U^\ast(x)=X\setminus Y$, hence $U^\ast$ has been extended to an ONA of $X$.
Let
$U$ be an ONA of $X$ finer than $U^\ast$ such that $X$ is $\kappa$-$D$ with
respect to $U$. We claim that $Y$ is $\kappa$-$D$ with
respect to $U\rest Y$. In this direction, let $D\subseteq Y$ of size $<\kappa$
be $U\rest Y$-sticky and let $y\in Y\setminus \bigcup U``D$. As $U$ is 
finer than
$U^\ast$, we have that $D$ is $U$-sticky, and hence there is $D'\ge_U D$ with 
$\card{D'}<\kappa$ and $y\in \bigcup U``D'$. Now $D''=D'\cap Y$ demonstrates that
$Y$ is $\kappa$-$D$.
$\eop_{\ref{closedkappaD}}$
\end{Proof of the Claim}

\begin{Note1} The choice of $U$ above depends on $Y$. This
and previous observations
motivate the following definition:
\end{Note1}

\begin{Definition}\label{uniformity} We say that $X$ is {\em uniformly}
$\kappa$-$D$ iff for every ONA $U^\ast$ of $X$, there is ONA $U$ finer than
$U^\ast$ such that every closed subspace of $X$ is $\kappa$-$D$ with respect to
$U$.
\end{Definition}

\rightline{$\between$}

The following argument is due to W. Fleissner.

\begin{Claim}\label{Bill} Suppose that $\kappa$ is a regular uncountable
cardinal and $X$ is a uniformly $\kappa$-$D$ space. Then $X$ is $D$.
\end{Claim}

\begin{Proof of the Claim} Let $U^\ast$ be an ONA of $X$ and let $U$ be a finer
ONA demonstrating that $X$ is uniformly $\kappa$-$D$. We shall show that $X$ is
$D$ with respect to $U$.

Given a $U$-sticky $D$ and $x\notin \bigcup U``D$. Let $Y\deq X\setminus\bigcup
U``D$,
so $Y$ is closed and $x\in Y$. As $\emptyset$ is $U$-sticky, we can find $D'$
with $\card{D'}<\kappa$ and $x\in \bigcup U``D'$, such that $D$ is $U\rest
Y$-sticky. Let $D''\deq D\cup D'$, then $D''\ge_U D$ is as required.
$\eop_{\ref{Bill}}$
\end{Proof of the Claim}

It is perhaps instructive to contrast Claim \ref{Bill} with the Conclusion of
Theorem \ref{upwards}.

\section{On downwards reflection} We prove several theorems which give conditions
on a $D$-space to have proper $D$-subspaces with specified properties. The
first theorem is an easy remark using the downward reflection of property
$D$ on closed subspaces, while the others are more involved and use a
large cardinal assumption.
As is often the case with such large cardinal
downward reflection arguments in topology, for one of the latter theorems an
additional assumption has to be made on the space in question in order to make
the reflection argument work. We concentrate on spaces with a small character, 
for a detailed discussion of other possible assumptions, the reader may consult
\cite{DoTaWe} and \cite{DoTaWe2}.

Let us first note that it is easy to obtain reflection 
results involving closed $D$-subspaces of a given $D$-space, because of the
following 

\begin{Observation}\label{closedxxx} If $X$ is a $D$-space and
$Y\subseteq X$ is a closed subspace of $X$, then $Y$ is a $D$-space.
\end{Observation}

\begin{Proof} The same proof as that of Claim \ref{closedkappaD}.
$\eop_{\ref{closedxxx}}$
\end{Proof}

\begin{Theorem}\label{closedsub} Suppose that $\kappa$ is a strong
limit and $X$ is a space of size $\kappa$.

{\noindent (1)} If $X$ is $D$, then for
every $\theta< \kappa$
and a subspace $Z$ of $X$ with $\card{Z}<\kappa$,
there is a closed $D$-subspace $Y$ of $X$ with $\theta<\card{Y}
<\kappa$ and $Z\subseteq Y$.

{\noindent (2)} If $X$ is $\sigma$-$D$ for some $\sigma<\kappa$, then for every
$\theta<\kappa$
and a subspace $Z$ of $X$ with $\card{Z}<\kappa$, there is a closed $\sigma$-$D$ subspace of $X$ with $\theta
<\card{Y}<\kappa$ and $Z\subseteq Y$.
\end{Theorem}

\begin{Proof} (1) Let $\lambda\in (\theta,\kappa)$ be a cardinal and let $Z$ be
a given subspace of $X$ of size
$<\kappa$, with an additional arbitrarily chosen set of $\lambda$ points of
$X$. Let $Y=\bar{Z}$. By Observation \ref{closedxxx}, it suffices to show that
$\card{Y}<\kappa$. This follows because $\card{Y}\le 2^{2^{\card{Z}}}<\kappa$,
as it is well known that for Hausdorff spaces $W$ we have $\card{\bar{W}}\le
2^{2^{\card{W}}}$, see \cite{Juhasz} II, 2.4.

{\noindent (2)} The same proof.
$\eop_{\ref{closedsub}}$
\end{Proof}

We thank W. Fleissner for simplifying the above argument and providing the gist
of the following

\begin{Theorem}\label{simplicity} (W. Fleissner, private communication) Suppose
that $\kappa$ is strongly inaccessible, and
that
$X$ is a locally $<\kappa$ $D$-space of size $\kappa$.

Then for every $\theta<\kappa$ and subspace $Z$ of $X$
with $\card{Z}<\kappa$, there is an open subspace $Y$ of $X$ with
$\theta<\card{Y}<\kappa$ which is $D$, and which contains $Z$ as a subspace.
\end{Theorem}

\begin{Proof} For any subspace $Z$ of $X$ with $\card{Z}<\kappa$, we shall
inductively choose sequences $\langle Y_n=Y_n(Z):\,n<\omega\rangle$
and $\langle Z_n=Z_n(Z):\,n<\omega\rangle$ as
follows:
\begin{itemize}
\item $Y_0=Z=Z_{0}$,
\item $Z_{n+1} = \bar{Y_{n}}$,
\item $Y_{n+1}$ is an open subspace of $X$ with $Y_{n+1}\supseteq
Z_{n+1}$ and $\card{Y_{n+1}}<\kappa$.
\end{itemize}
Let $Y^\ast=Y^\ast(Z)=\bigcup_{n<\omega}Y_n=\bigcup_{n<\omega}Z_n$.
It suffices to show that $Y^\ast(Z)$
can always be chosen as required, and that it is an open $D$ subspace of $X$ of
size $<\kappa$. Let $Z$ with $\card{Z}<\kappa$ be given.

It is clear that the sequence of $Z_n$s can be chosen as required, so
we show by induction on $n<\omega$ that we can choose $Y_n$s as
well. A part of the inductive hypothesis is that $\card{Y_n}<\kappa$.
Coming to $Y_{n+1}$, we have $\card{Z_{n+1}}<2^{2^{\card{Y_{n}}}}
<\kappa$, as
$\kappa$ is a strong limit. Then to choose $Y_{n+1}$, we pick for every
$y\in Z_{n+1}$ an open neighbourhood of size $<\kappa$ and let
$Y_{n+1}$
be the union of all these. As $\kappa$ is regular, we have
$\card{Y_{n+1}} <\kappa$.

It is clear that $\card{Y^\ast}<\kappa$ and that $Y^\ast$ is open. We
show that $Y^\ast$ is $D$, so let $U$ be a given ONA of
$Y^\ast$.

For $y\in Y^\ast$ let $n(y)$ be the first $n$ such that $y\in Z_n$ and let
\[
V(y)=U(y)\setminus \bigcup_{k<n(y)} Z_k.
\]

Hence $V(y)$ is open and $V(y)\subseteq U(y)$, while
$y\in V(y)$. We shall find a $V$-sticky $D^\ast$
such that $\bigcup V``D^\ast=Y^\ast$, which is clearly sufficient. Note that 
by the choice of $V$ we have that if $n<k$ and $D\subseteq Z_n$ is $V\rest
Z_n$-sticky, then $D$ is $V\rest Z_k$-sticky, and also $D$ is $V$-sticky.

Now we build by recursion on $n<\omega$ an increasing sequence $\langle
D_n:\,n<\omega\rangle$ of elements of ${\Bbf P}_V$
such that each $D_n$ is $V\rest Z_n$-sticky and $\bigcup
V``D_n\supseteq Z_n$. This can be done because each $Z_n$ is closed, using Theorem
\ref{FlSt}. At the end let $D^\ast=\bigcup_{n<\omega}D_n$.
$\eop_{\ref{simplicity}}$
\end{Proof}

This argument greatly simplifies the one we had originally, which used a
measurable cardinal. Note that the conclusion of the theorem clearly implies
that the assumption of $X$ being locally $<\kappa$ is necessary.
The same method works for uniformly $\sigma$-$D$-spaces, but does not seem to 
work if uniformity is not assumed.
In order to obtain somewhat of an an analogue for $\sigma$-$D$-spaces, we need a
stronger large
cardinal assumption, and an additional  assumption on $X$.

\begin{Theorem}\label{downaward} Suppose that $\kappa$ is
measurable and $X$ is a locally $<\kappa$ space of size $\kappa$, such that
every point in $X$ has a point-base of size $<\kappa$, and suppose that
$X$ is $\sigma$-$D$ for some $\sigma<\kappa$.

Then for every $\theta<\kappa$ there is $Y\subseteq X$ which is an open
$D$-subspace of $X$ and satisfies
\[
\theta<\card{Y}<\kappa.
\]
\end{Theorem}

\begin{Proof} Let $j:\,{\bf V}\into M$ be an embedding witnessing
that $\kappa$ is measurable, so in particular $\kappa$ is the critical point
of $j$ and ${}^\kappa M\subseteq M$.

Suppose that $(X,\tau)$ is a given space with the properties as listed above.
The idea of the proof is, as one would imagine, that $j``X$ is a subspace of
$j(X)$ that has the properties as required of $Y$ when translated by $j$, so
that by elementarity $X$ must have a subspace $Y$ as required. 
However,
topological reflections arguments are not so simple
as many notions involved
are highly non-first order, and in
particular, an assumption has to be used to guarantee that $j``X$ is
actually a subspace of $j(X)$. In fact, we first have to clarify which
topologies we have in mind when discussing $j(X)$ and $j``X$.

For simplicity we shall assume,
without loss of generality, that $X$ as a set is $\kappa$. 
By the assumption on the character of $X$, we can fix a
sequence $\bar{B}=\langle \bar{B}_\alpha:\,\alpha<\kappa\rangle$
such that for every $\alpha<\kappa$ we have that
$\bar{B}_\alpha=\langle B^\alpha_\zeta:\,\zeta<\zeta_\alpha<\kappa\rangle$
is a point base for $\alpha\in X$. As $X$ is locally $<\kappa$, we may assume
that $\card{B^\alpha_\zeta}<\kappa$ for all $\alpha$ and $\zeta$.
Let  \[
\BB=\{B^\alpha_\zeta:\,\alpha<\kappa,
\zeta<\zeta_\alpha\},
\]
so $\BB$ is a basis for the topology $\tau$ on $X$.
By elementarity, $j(\BB)$ is a basis for a topological space in $M$
whose set of points is $j(\kappa)$.
We shall abbreviate this space as $j(X)$.
For every
$\alpha<\kappa$, a point-base at $j(\alpha)=\alpha$ in $j(X)$ is
given by $\langle j(B^\alpha_\zeta):\,\zeta<\zeta_\alpha\rangle$. 
As for every $\zeta<\zeta_\alpha$ we have that
$j(B^\alpha_\zeta)\cap j``X(=\kappa)=j``B^\alpha_\zeta=B^\alpha_\zeta$, we have
that $\BB=j``\BB$ generates the subspace topology on $j``X$. Clearly the
original topology of $X$ is at least as fine as this topology, but in fact 
our assumptions guarantee that these two topologies are the same. For if
$Y\subseteq X$ is open, then it is the union of a sequence of $\le\kappa$
elements of $\BB$, hence this sequence is a member of $M$ and so $Y$ is open
in $M$. 

\begin{Observation}\label{openness} In $M$, $X$ is an open subspace
of $j(X)$.
\end{Observation}

\begin{Proof of the Observation} Let $x\in X$, and let $\zeta<\zeta_x$
be arbitrary. Then $\card{B^x_\zeta}<\kappa$ and $x\in j(B^x_\zeta)=
j``B^x_\zeta=B^x_\zeta\subseteq X$, hence $X$ is open.
$\eop_{\ref{openness}}$
\end{Proof of the Observation}

\begin{Observation}\label{Dsubspace}  In $M$ we have that
$X$ is a $\sigma$-$D$-subspace of $j(X)$.
\end{Observation}

\begin{Proof of the Observation} In $M$, let $U$ be an ONA of $X$,
Back in ${\bf V}$, as $X$ is $\sigma$-$D$, there is a finer ONA $W$ of $X$ such
that $X$ is $\sigma$-$D$ with respect to $W$. Remembering that the topologies
of $X$ in $M$ and ${\bf V}$ are the same and that ${}^\kappa M\subseteq M$,
we have that $W\in M$ and clearly $W$ is a finer ONA of $X$ then $U$ is. We
claim that $M$ satisfies that $X$ is $\sigma$-$D$ with respect to $W$. 

We work in $M$ and let $x\in X$ and $D\subseteq X$ be such that
$\card{D}<\sigma$ and $D$ is $W$-sticky. Then in ${\bf V}$ we
have that $D$ is $W$-sticky, so there is $D'\ge D$ with
$x\in D'$ and $\card{D'}<\sigma$. It is easy to verify that $D'\in M$
satisfies the same requirements.
$\eop_{\ref{Dsubspace}}$
\end{Proof of the Observation}

Now let us finish the proof of the Theorem. 
Let $\theta<\kappa$ be given. As
$X\in M$, in $M$ we have
that there is $Y=X\subseteq j(X)$ of size 
a cardinal in $(j(\theta)=\theta,j(\kappa))$,
such that for every $U$ which is a function from $Y$ to
$\{O\cap Y:\,O\in j(\BB)\}$ with the property $y\in U(y)$ for all $y$,
there is a function $W$ from $Y$ to $\{O\cap Y:\,O\in j(\BB)\}$ with the
following properties: 

\begin{description}
\item{(i)} $(\forall y\in Y)\,[W(y)\subseteq U(y)\,\,\&\,\,y\in W(y)]$,

\item{(ii)} 
Let $\varphi(D;Y, j(\BB), W)$ stand for 
\[
\card{D}<\sigma\,\,\&\,\,(\exists\FF\subseteq \{O\cap Y:\,O\in
j(\BB)\})(Y\setminus D =\bigcup\FF)\,\,\&\,\,
\]
\[
(\forall y\in D)(\exists B\in j(\BB))
(B\cap D=\{y\})\,\,\&\,\,
\]
\[
(\forall y\in Y)[W(y)\cap D\neq\emptyset\implies
(\exists z\in D)y\in W(z)].
\]
Then 
\[
(\forall x\in Y)(\forall D\subseteq
Y)[\varphi(D;Y, j(\BB), W)\implies (\exists D'\supseteq D)
\]
\[(x\in D'\,\,\&\,\,
\varphi(D';Y, j(\BB), W)\,\,\&\,\,(\forall z\in D'\setminus D)(W(z)\cap
D=\emptyset).
\]
\end{description}

By elementarity, in ${\bf V}$ we can find a subset $Y$ of $X$ with
$\theta<\card{Y}<\kappa$
such that for every $U$ which is a function from $Y$ to
$\{O\cap Y:\,O\in \BB\}$ with the property $y\in U(y)$ for all $y$,
there is $W$ from $Y$ to $\{O\cap Y:\,O\in \BB\}$

with the following properties:
\begin{description}
\item{(i)} as above
\item{(ii)} \[
(\forall x\in Y)(\forall D\subseteq
Y)[\varphi(D;Y, \BB, W)\implies
\]
\[
(\exists D'\supseteq D)(x\in D'\,\,\&\,\,
\varphi(D';Y, \BB, W)\,\,\&\,\,(\forall z\in D'\setminus D)(W(z)\cap
D=\emptyset),
\]
\end{description}
which is as required.
$\eop_{\ref{downaward}}$
\end{Proof}

At this point it is natural to ask if we can obtain similar downward transfer
properties between $\aleph_2$ and $\aleph_1$ by applying the technique of
generic embeddings, for example by a L\'evy collapse of a large cardinal to
$\aleph_2$, or using a huge cardinal embedding in the fashion of
M. Foreman and R. Laver in \cite{FoLa}. Although some weak partial results can
be easily obtained, difficulties with the transfer of property $D$ make it
unclear whether the exact analogue of any of our downward reflection
principles can be
obtained. Let us also note that in Theorem \ref{downaward} one could perhaps relax the
various assumptions made, but we feel that the version presented is convenient
as a contrast to the theorems from Section \S1, and for the simplicity of
reflection arguments, while the details of the consistency strength and
topological strength investigation might be premature before we understand more
about the relationship between $D$, $\sigma$-$D$ and uniformly $\sigma$-$D$.

\section{Discrete families of sets} The $D$-space problem can be formulated
as a purely combinatorial statement involving discrete families
of sets, as will be
shown below, where we shall also exhibit some basic properties of the families
in question. We commence by a definition.

\begin{Definition} (1) A non-empty family $\FF$ of non-empty sets is said to
be {\em discrete} iff there is a choice function $f$ on $\FF$ such that
\[
F_0\neq F_1\in \FF\implies f(F_0)\notin F_1.
\]
A function $f$ as above is called a {\em discretisation} of $\FF$.

{\noindent (2)} Let $\FF$ be as above, and $\GG\subseteq\PP(\bigcup \FF)$. We say
that $\FF$ is $\GG$-discrete iff there is a ${\cal D}\subseteq \FF$
and a discretisation $f$ of ${\cal D}$ such that
\begin{description}
\item{(i)} $\bigcup{\cal D}=\bigcup\FF$
\item{(ii)} $\{f(D):\,D\in{\cal D}\}\in \GG$.
\end{description}
In such a case, the pair $({\cal D},f)$ is called a $\GG$-{\em discretisation}
of $\FF$.
\end{Definition}
\rightline{$\between$}
Hence a (Hausdorff) topological space $X$ is $D$ iff every ONA $U$ of $X$, 
there is a closed discretisation of $\{U(x):\,x\in X\}$. The definition
of discretisation is similar to the definition of a transversal,
which is a one-to-one choice function. However, one
should note that the requirement
for a function to be a discretisation is stronger than just being
one-to-one. Of course, the two notions coincide when the family $\FF$ consists
of pairwise disjoint elements, but we are mainly interested in the
cases where the existence of a transversal is not an obvious consequence of the
axiom of choice. There is a body of work about
transversals, cf. Shelah's book \cite{Sh g}, often concentrating on the 
incompactness properties, that is,
families in which every smaller subfamily has a transversal, but the whole
family does not have it. Such problems are known to be equivalent to the
existence of certain families of functions, see II 6.2 in \cite{Sh g}.
A similar argument can be used to characterise the existence of a
non-discrete family of sets in which every smaller subfamily is discrete,
as expressed by the following Theorem \ref{npd}.
Although the proof is very much the same as that of the corresponding one
in the case of transversals, in \cite{Sh g}, as the details there are
not fully explained and
as we need them for later use, we have decided to
spell out the proof here. Following this theorem we shall obtain
as a corollary a connection between families whose small subfamilies
have transversals and such families that in addition satisfy a disreteness
requirement (see Theorem \ref{novi}). Let
us first make

\begin{Observation}\label{never} If $\lambda$ is a cardinal, there is no family
of $>\lambda$ subsets of $\lambda$ that has a transversal.
\end{Observation}

\begin{Proof} Suppose that ${\cal A}=\{A_\alpha:\,\alpha<\alpha^\ast\}$
is such a family with $\card{\alpha^\ast}>\lambda$. Let $f$ be a transversal
of ${\cal A}$. Then $\{f(A_\alpha):\,\alpha<\alpha^\ast\}$ is a subset of
$\lambda$ of size $>\lambda$, a contradiction.
$\eop_{\ref{never}}$
\end{Proof}

\begin{Theorem}\label{npd} Suppose that $\mu>\lambda\ge\theta\ge\kappa$ are
infinite cardinals.
Then the following are equivalent:

(A) There is a family $\PP^\ast$ of $\mu$ subsets of $\lambda$, each of power
$\le \kappa$ 
(none of whose subfamilies of size $\mu$ has a transversal, but) whose every
subfamily of size $<\theta$ is discrete,

and

(B) There is a regular ideal $J$ on $\kappa$, and a family $\FF^\ast$ of $\mu$
many functions from $\kappa$ to $\lambda$ such that for every subfamily $\FF$
of $\FF^\ast$ with $\card{\FF}<\theta$, there is a sequence
$\langle s_f:\,f\in \FF\rangle$ of sets in $J$ such that 
\[
i\in \kappa\setminus (s_f\cap s_g)\implies f(i)\neq g(i),
\]
but there is no such sequence for any subfamily of $\FF^\ast$ which has
size $\mu$. Moreover, if $\FF\subseteq\FF^\ast$ is of size $\ge
\lambda^+$ there is no
sequence $\langle s_f:\,f\in \FF\rangle$ of sets in $J$ such that
\[
i\in \kappa\setminus (s_f\cup s_g)\implies f(i)\neq g(i).
\]
\end{Theorem}

\begin{Note1} 
Claim II 6.2 of \cite{Sh g} gives a similar characterisation, in which
the existence of a discretisation is replaced by the existence of a
transversal, and ``$s_f\cap s_g$" in (B) above is replaced by ``$s_f\cup s_g$".
\end{Note1}

\begin{Proof} (A)$\implies$(B). Let us enumerate $[\kappa]^{<\aleph_0}$
as $\{w_i:\,i<\kappa\}$, and let $F$ be a bijection from
$[\lambda]^{<\aleph_0}$ onto $\lambda$. We define
\[
J\deq\{A\subseteq\kappa:\,(\exists i<\kappa)(\forall j\in A)w_i\nsubseteq
w_j\}.
\]
It is clear that $J$ is a proper ideal on $\kappa$. If $A\subseteq \kappa$ is
bounded, with $\sup(A)=\alpha<\kappa$, then $\card{\cup\{w_j:\,j\le\alpha\}}
<\kappa$, so there must be an $i<\kappa$ such that $w_i$ is not a subset of
$w_j$ for any $j\in A$. Hence, $J$ is regular.

Let $\PP^\ast$ be a family as in the assumptions of (A). For $X\in \PP^\ast$,
let us enumerate $X=\{\alpha^X_\zeta:\,\zeta<\zeta_X\le\kappa\}$. For $i<\kappa$,
let 
\[
f_X(i)\deq F(\{\alpha^X_j:\,j\in w_i\cap\zeta_X\}),
\]
hence each $f_X$ is a function from $\kappa$ into $\lambda$. Let $\FF^\ast
\deq\{f_X:\,X\in \PP^\ast\}$, and let us claim that $\FF^\ast$ is as required.
First note that $X_0\neq X_1\implies f_{X_0}\neq f_{X_1}$, so the size
of $\FF^\ast$ is $\mu$.

Let $\FF\subseteq \FF^\ast$ be of size $<\theta$, so $\FF=\{f_X:\,X\in \PP\}$
for some $\PP\subseteq\PP^\ast$ with $\card{\PP}<\theta$. Hence $\PP$ is
discrete, and we can fix a discretisation $h$ of $\PP$. In particular
notice that $h(X)\in X$ for all $X\in \PP$. We define
\[
s_{f_X}\deq\{i<\kappa:\,h(X)\notin \{\alpha^X_j:\,j\in w_i\}\},
\]
for $X\in \PP$. Notice that each $s_{f_X}\in J$, as one can take $i<\kappa$ such that
$h(X)=\alpha^X_{\zeta^\ast}$ for some
$\zeta^\ast$ and $w_i=\{\zeta^\ast\}$. Then if $j\in
s_{f_X}$, we have that $\{h(X)\}$
is not contained in $\{\alpha^X_\zeta:\,\zeta\in w_j\}$,
so $\zeta^\ast\notin w_j$, and hence $w_i\nsubseteq w_j$.

If $X\neq Y\in \PP$, then we have that $h(X)\notin Y$, so
clearly $h(X)\notin \{\alpha^Y_j:\,j\in w_i\}$ for any $i$. But if $i\notin
s_{f_X}$, we have $h(X)\in \{\alpha^X_j:\,j\in w_i\}$, so
\[
\{\alpha^X_j:\,j\in w_i\cap \zeta_X\}\neq \{\alpha^Y_j:\,j\in w_i\cap\zeta_Y\},
\]
and in particular the images of these sets under $F$ are distinct. Hence, for
each such $i$ we have $f_X(i)\neq f_Y(i)$, and by symmetry
the same is true for $i\notin s_{f_Y}$.

Now let us prove the last claim of (B).
Suppose that $\FF_0\subseteq \FF^\ast$ has size $\ge\lambda^+$ and that $\langle
s_{f_X}:\,X\in \FF_0\rangle$ can be defined as required. For each
$X\in \FF_0$, we can find $i_X\in (\kappa\setminus s_{f_X})$, 
which is possible as $J$ is proper. Then there is
$\FF_1\subseteq \FF_0$ of size $\ge\lambda^+$ and $i^\ast<\kappa$ such that for all
$X\in \FF_1$ we have $i_X=i^\ast$. But then $f_X(i^\ast)$ for $X\in \FF_1$
are $\ge\lambda^+$ distinct elements of $\lambda$, a contradiction.
\smallskip

{\noindent (B)$\implies$(A)}. This direction is easier: fix a bijection
$F$ between $\kappa\times\lambda$ and $\lambda$. Starting with
$\FF^\ast$ as in (B), for $f\in \FF^\ast$ define
$X_f=\{F((i,f(i)):\,i<\kappa\}$. Let $\PP^\ast\deq\{X_f:\,f\in \FF^\ast\}$,
it is easy to check that $\PP^\ast$ is as required.

$\eop_{\ref{npd}}$
\end{Proof}

\begin{Note1} A non-trivial $\Delta$-system is an example of a discrete
family. Hence, for the Theorem \ref{npd} to be interesting, we need at least to
be in a situation in which $\Delta$-system Lemma between $\mu$ and $\kappa^+$
does not hold, so we should have $\sigma^\kappa\ge\mu$ for some $\sigma<\mu$.
\end{Note1}

The following theorem establishes a simple direct relationship between
transversals and discretisations.

\begin{Theorem}\label{novi} Suppose that $\mu>\lambda\ge\theta>\kappa$ are
infinite cardinals, and suppose that $\PP$ is a family of $\mu$ elements
of $[\lambda]^{\le\kappa}$ such that every subfamily of $\PP$ of size
$<\theta$ has a transversal. Then
there is such a family $\PP$ such that in addition, every subfamily of $\PP$
whose size $\sigma$ satisfies

$\sigma<\theta$ and $\cf(\sigma)>\kappa$, has a subfamily of size $\cf(\sigma)$
that has a discretisation. 
\end{Theorem}

\begin{Proof} Let $\PP$ satisfy the assumptions of the theorem.
Then Shelah's result from \cite{Sh g} II 6.2 is that the proof of
Theorem \ref{npd} with the ideal $J$ defined as there
but using $\PP$ in place of $\PP^\ast$, yields a family
$\FF^\ast$ of $\mu$ many functions from $\kappa$ to $\lambda$
such that for every subfamily $\FF$ of $\FF^\ast$ of size $<\theta$
there is a sequence $\langle s_f:\,f\in \FF\rangle$ of elements of $J$ such that
\begin{equation*}
\card{\{g\in \FF:\,(\exists i \in \kappa\setminus (s_f\cup s_g))
f(i)=g(i)\}}\le\kappa.\tag{$\ast$} \end{equation*}
We shall show that $\FF^\ast$ has the property that for every subfamily $\FF$
of $\FF^\ast$ of size $\sigma<\theta$ with $\cf(\sigma)>\kappa$ there is a
subfamily $\FF'$ of $\FF$ of size $\cf(\sigma)$ for which $(\ast)$ holds with
$``s_f\cup s_g"$ replaced by $``s_f\cap s_g$" and $``\le\kappa"$ replaced
by ``$\le 1$".

Suppose that $\FF$ is a subfamily of $\FF^\ast$ of size $\sigma<\theta$ with 
$\cf(\sigma)>\kappa$. Let $\langle s_f:\,f\in \FF\rangle$ be
as guaranteed by $(\ast)$. Note that the sets
\[
A_i\deq\{j<\kappa:\,w_i\nsubseteq w_j\}
\]
generate the ideal $J$, in the sense that every element of $J$ is a subset of
an $A_i$. Here the sets $w_i$ are as
in the proof of Theorem \ref{npd}.
Since the conclusion of $(\ast)$ does not change if we replace each
$s_f$ by a set that is larger than $s_f$ but still in $J$, we can assume that
for each $f$ there is $i(f)<\kappa$ such that $s_f=A_{i(f)}$. Hence there is
$i^\ast<\kappa$ such that for $\cf(\sigma)$ many $f$ we have $i(f)=i^\ast$.
For such $f$ let $s^\ast=s_f$. Let ${\FF}_0=\{f\in \FF:\,s_f=s^\ast\}$, hence
$\card{{\FF}_0}=\cf(\sigma)$ and we can enumerate
$\FF_0=\{f_\zeta:\,\zeta<\cf(\sigma)\}$.

Given $\zeta<\cf(\sigma)$, we notice that 
\[
d(\zeta)\deq\sup\{\xi<\cf(\sigma):\,(\exists i\in
\kappa\setminus s^\ast) f_\xi(i)=f_\zeta(i)\}<\cf(\sigma).
\]
Hence, by induction on $\alpha<\cf(\sigma)$ we can define $\zeta_\alpha$ as
follows:

Let $\zeta_0=0$. Given $\zeta_\alpha$ let $\zeta_{\alpha+1}=d(\zeta_\alpha)+
1$.
For $\alpha$ a limit ordinal $<\cf(\sigma)$ let $\zeta_\alpha\deq
\sup\{\zeta_\beta:\,\beta<\alpha\}$.
Then we can let $\FF'=\{f_{\zeta_\alpha}:\,\alpha<\cf(\sigma)\}$.

Having shown this property of $\FF^\ast$, we proceed to define
$\PP^{\ast}$ as in the proof of $(B)\implies (A)$ of Theorem \ref{npd}.
It is easy to check that $\PP^\ast$ is as required.
$\eop_{\ref{novi}}$
\end{Proof}

Shelah uses the following notation for the situation in the assumptions
of the Theorem \ref{novi}

\begin{Definition} Suppose that $\mu\ge\lambda\ge\theta_1\ge\theta_2+\kappa$
and $J$ is an ideal on $\kappa$. Then ${\rm
NPT}_J(\mu,\lambda,\theta_1,\theta_2,\kappa)$ means that there is a family
$\FF^\ast$ of $\mu$ functions from $\kappa$ to $\lambda$ such that 
\begin{description}
\item{(a)} for any subfamily $\FF$ of $\FF^\ast$ with $\card{\FF}<\theta_1$,
there is a sequence $\langle s_f:\,f\in \FF\rangle$ of elements of $J$ such
that for each $f\in \FF$
\[
\card{\{g\in \FF:\,(\exists i\in \kappa\setminus (s_f\cup s_g))f(i)=g(i)\}}
<\theta_2,
\]
\item{(b)} the analogue of (a) with $\FF^\ast$ in place of $\FF$ fails.
\end{description}
\end{Definition}
\rightline{$\between$}

In conjunction with the following Theorem \ref{when} of Shelah from
\cite{Sh g} II 6.3, Theorem \ref{novi}
can be used to read off non-trivial instances of families with discretisations.
The function cov is discussed in detail in \cite{Sh g}, but an instance of 
a situation in which the assumptions of Theorem \ref{when} hold is
\[
\lambda^{\aleph_0}>\lambda^+\,\,\&\,\,(\forall\theta<\lambda)(\theta^{\aleph_0}
<\lambda).
\]

\begin{Theorem}\label{when} (Shelah) Suppose that $\lambda>\cf(\lambda)=\aleph_0$
and ${\rm cov}(\lambda,\lambda,\aleph_1,2)>\lambda^+$.
Then ${\rm NPT}_{{\rm Fin}}(\lambda^+,\lambda,\lambda^+,2,\aleph_0)$
holds.
\end{Theorem}

In the above Fin stands for the ideal of finite subsets of $\omega$.

\section{Concluding Remarks} We investigated reflection phenomena that arise in 
connection with van Douwen's notion of $D$-spaces. In the first two sections
we concentrated on the topological aspects of this problem, studying both 
upwards and downwards reflection. The last section shows that there is a
purely combinatorial aspect of the problem, in the sense that one can define
a generalisation of $D$-property that is formulated in terms of a covering
of one family of sets by another. Then one can talk about discreteness
properties of such covers and obtain the original topological formulation
of $D$-spaces as a particular instance of this more general setting.
In tune with the rest of the paper, we concentrated again on reflection 
properties of such covers and showed that such properties of discrete families of sets
have a strong connection with the well studied 
combinatorial problem of the existence of transversals. 
This indicates that it would be of interest to study discrete families of sets from
the purely combinatorial point of view, an investigation that is outside of the
scope of this paper.
One could then however hope that families of sets with
given discreteness properties could be topologised in order to give examples 
of topological spaces 
of some relevance to the $D$-space problem.

\eject


\begin{thebibliography}{99}

\bibitem{Lutzer}  E.~K. van Douwen and D.~J. Lutzer, {\em A note on the 
paracompactness in generalized linearly ordered spaces\/}, Proceedings of the
American Mathematical Society, {\bf 125} No.4 (1997), 1237--1245.

\bibitem{DoPf} E.K. van Douwen and W.F. Pfeffer, {\em Some properties
of the Sorgenfrey line and related spaces\/}, Pacific Journal of Mathematics,
{\bf 81} No.2 (1979), 371--376.

\bibitem{DoTaWe} A.Dow, F.D. Tall and W.A.R. Weiss, {\em New proofs
of the consistency of the normal Moore space conjecture I\/},
Topology and its Applications, {\bf 37} (1990), 3--51.

\bibitem{DoTaWe2} A.Dow, F.D. Tall and W.A.R. Weiss, {\em New proofs
of the consistency of the normal Moore space conjecture II\/},
Topology and its Applications, {\bf 37} (1990), 115--129.

\bibitem{Adrienne} W. Fleissner and A.~M. Stanley, {\em $D$-spaces\/},
Topology and its Applications, {\bf 114} No.3 (2001), 261--271.

\bibitem{FoLa} M. Foreman and R. Laver, {\em Some Downwards Transfer
Properties for $\aleph_2$\/}, Advances in Mathematics, {\bf 67} No. 2
(1988), 230--238.

\bibitem{Juhasz} I. Juh\'asz, {\em Cardinal Functions in Topology-Ten Years
Later\/}, Mathematical Centre Tracts, {\bf 123}, Amsterdam 1980, 160+iv pp.

\bibitem{JuShSo} I. Juh\'asz, S. Shelah and L. Soukup, {\em More on 
countably compact, locally countable spaces\/}, Israel Journal of Mathematics,
{\bf 62}, No. 3 (1988), 302--310.

\bibitem{Sh g} S. Shelah, {\em Cardinal Arithmetic\/}, Oxford University Press,
 Oxford 1994, 481 + xxxi pp.

\bibitem{Sh80} S. Shelah, {\em A weak generalization of MA to higher
cardinals\/}, Israel Journal of Mathematics, {\bf 30} No. 4 (1978), 297--306.

\bibitem{counter} L.A. Steen and J.A. Seebach Jr., {\em Counterexamples
in Topology\/}, Dover edition 1995 (1st edition 1970), 244 + xi pp.





\end{thebibliography}
\end{document}